\documentclass[11pt]{amsart}

\usepackage{amsmath}
\usepackage{amssymb}
\usepackage{graphicx}

\usepackage{amssymb,mathrsfs,amstext, enumerate, comment}
\usepackage[noadjust]{cite}
\usepackage[plainpages=false,colorlinks,hyperindex,pdfpagemode=None,bookmarksopen,linkcolor=blue,citecolor=blue,urlcolor=blue]{hyperref}

\usepackage{multicol}

\newtheorem{thm}{Theorem}[section]

\newtheorem{prop}[thm]{Proposition}
\newtheorem{cor}[thm]{Corollary}

\theoremstyle{definition}
\newtheorem{defn}[thm]{Definition}

\numberwithin{equation}{section}

\def\M{{\mathscr{ M}}}

\def\A{{\mathscr{A} }}

\def\M{{\mathscr{M} }}

\newcommand{\cA}{\mathcal{A}}

\newcommand{\cF}{\mathcal{F}}
\newcommand{\cB}{\mathcal{B}}

\newcommand{\bM}{\mathbf{M}}
\newcommand{\cS}{\mathcal{S}}

\begin{document}

\title[Operational $2$-local automorphisms/derivations]{Operational  $2$-local automorphisms/derivations}

\author[L. Wang]{Liguang Wang}
\address[L. Wang]{School of Mathematical Sciences, Qufu normal University, Qufu 273165, China.}
\email{\tt wangliguang0510@163.com}

\author[N.-C. Wong]{Ngai-Ching Wong}
\address[N.-C. Wong]{Department of Applied Mathematics, National Sun Yat-sen University, Kaohsiung, 80424, Taiwan, \and
	School of Mathematical Sciences, Tiangong University, Tianjin 300387, China.}
\email{\tt wong@math.nsysu.edu.tw}



\thanks{Corresponding author: Ngai-Ching Wong, wong@math.nsysu.edu.tw}

\begin{abstract}
Let $\phi: \cA\to \cA$ be a (not necessarily linear, additive or continuous)
 map of a standard operator algebra.
Suppose for any $a,b\in \cA$ there is an algebra automorphism $\theta_{a,b}$  of $\cA$ such that
\begin{align*}
\phi(a)\phi(b) = \theta_{a,b}(ab).
\end{align*}
We   show that either $\phi$ or $-\phi$   is a linear Jordan homomorphism.
 Similar results are obtained when any of the following conditions is satisfied:
\begin{align*}
\phi(a) + \phi(b) &= \theta_{a,b}(a+b), \\
\phi(a)\phi(b)+\phi(b)\phi(a) &= \theta_{a,b}(ab+ba), \quad\text{or} \\
\phi(a)\phi(b)\phi(a) &= \theta_{a,b}(aba).
\end{align*}

We also show that
a  map $\phi: \M\to \M$ of a semi-finite
von Neumann algebra $\M$ is a linear derivation if
for every $a,b\in \M$ there is a linear
derivation $D_{a,b}$ of $\M$ such that
$$
\phi(a)b + a\phi(b) = D_{a,b}(ab).
$$
\end{abstract}

 \maketitle

\section{Introduction}

Let $\cA$ be a complex ($*$-)algebra.
We call a  map
  $\theta: \cA\to \cA$ a  ($*$-) \emph{automorphism} if $\theta$ is bijective, ($*$-)linear,  and multiplicative.
We call $\theta$ a \emph{local} ($*$-)\emph{automorphism} if for every $a$ in $\cA$ there is a ($*$-)automorphism
$\theta_a$ of $\cA$, depending on $a$, such that
$
\theta(a)=\theta_a(a).
$
Although a local automorphism preserves
idempotents, square zero elements, central elements, invertibility, and
 spectrum, it can be nonlinear and/or non-multiplicative (see, e.g., \cite{LW07}).

 The notion of local
automorphisms is introduced by Larson and Sourour \cite{LS1990}.
They showed that every surjective linear local automorphism of a
matrix algebra is either an automorphism or an anti-automorphism,
and that of $B(H)$ is an automorphism whenever $H$ is an infinite
dimensional Hilbert space (see also Bre\v{s}ar and \v{S}emrl \cite{bresar95}).
It is also known that a surjective linear local automorphism  of a (resp.\ purely infinite)  $C^*$-algebra of real rank zero,
or a (resp.\ properly infinite)  von Neumann algebra, is a linear
Jordan isomorphism (resp.\ an automorphism), while
a linear local automorphism  of an abelian C*-algebra is always an algebra homomorphism (see, e.g., \cite{LW07}).

We call a  (not necessarily linear, additive or continuous) map
 $\theta:\cA\to \cA$ a \emph{$2$-local} ($*$-)\emph{automor\-phism} if for any  $a,b$ in $\cA$ there is a ($*$-)automorphism
$\theta_{a,b}$ of $\cA$ such that $\theta(a)=\theta_{a,b}(a)$ and  $\theta(b)=\theta_{a,b}(b)$.
It is clear that a
 $2$-local ($*$-)automorphism $\theta$ is a local ($*$-)automorphism.
It is also known that every $2$-local automorphism
 of a standard operator algebra $\cA$ on a Banach space
is an algebra homomorphism provided
that  it is continuous with respect to the weak operator topology or its range contains all  finite rank operators,
while a surjective $2$-local automorphism
of a C*-algebra is an algebra automorphism (see \cite{LW06}).
Moreover, every $2$-local $*$-automorphism of a von Neumann
algebra  is an algebra $*$-homomorphism (\cite{BPGP15}).

On the other hand,
a linear map $D: \A\to \A$ is said to be a \emph{derivation} if
$$
D(ab) = D(a)b + aD(b) \quad\text{for all $a,b\in \A$.}
$$
Every derivation on a C$^{\ast}$-algebra is (norm) continuous (see \cite[Lemma 4.1.3]{Saki71}).
Kadison and Sakai showed that every  (linear) derivation $D: \M\to \M$ of a von Neumann algebra $\M$ is inner, namely, there is
an $s\in \M$ such that
$$
D(a) = sa - as \quad\text{for all $a\in \M$}
$$
 (see \cite[Theorem 4.1.6]{Saki71}).
It is also known that every derivation of a simple $C^{\ast}$-algebra  with identity is inner (\cite{Saki68}; see also
 \cite[Theorem 4.1.11 and Corollary 4.1.7]{Saki71}),
 and  every derivation from a standard
operator algebra into $\cB(X)$ is an inner derivation, where $X$ is a Banach space (\cite{Chernoff73}).

We call a  map
 $\phi: \A\to \A$   a \emph{local derivation} if
for any $a$ in $\A$ there is a  derivation $D_a$ of $\A$, depending on $a$, such that
$
\phi(a) = D_a(a). 
$
Every linear local derivation of a $C^\ast$-algebra is continuous (\cite[Theorem 7.5]{Johnson01}), and
indeed a derivation (\cite[Corollary 1]{shulman94}; see also \cite{bresar92}).
In particular, every  {linear} local derivation $\phi$ on a von Neumann algebra  is  an inner derivation (\cite[Theorem A]{Kad90}).

We call a  map $\phi: \A\to \A$ a   \emph{$2$-local derivation} if
for any $a,b$ in $\A$ there is  a derivation $D_{a,b}$ of $\cA$ such that
$\phi(a)=D_{a,b}(a)$  and $\phi(b)=D_{a,b}(b)$.
It is clear that $2$-local derivations are local derivations.
Ayupov and Kudaybergenov showed that
every $2$-local derivation of a von Neumann algebra is
a linear local derivation, and thus
an inner derivation (\cite[Theorem 2.1]{AK15}).

Recently, Moln\'{a}r \cite{Molnar22} studied some operational forms of $2$-local automorphisms/ deri\-vations.
They are those maps $\phi: \A\to \A$ satisfying
one of the following conditions: for any $a,b\in \A$, there
is an automorphism $\theta_{a,b}$ or a (linear) derivation $D_{a,b}$ of $\cA$ such that
\begin{align}
\phi(a) + \phi(b) &= \theta_{a,b}(a+b), \label{eq:w2-local+}\\
\phi(a)\phi(b) &= \theta_{a,b}(ab), \label{eq:w2-localx} \quad\text{or}\\
  \phi(a)b + a\phi(b) &= D_{a,b}(ab) \label{eq:w2-localD}.
\end{align}
It is clear that any $2$-local automorphism/derivation  is an operational $2$-local automorphism/derivation.

Although Proposition \ref{prop:Molnar}(c) below is stated for standard operator algebras on  Hilbert spaces
 in \cite[Theorem 2.7]{Molnar22}, the
statement and the proof there are indeed also good for the case when $E$ is a Banach space.

\begin{prop}[Moln\'{a}r {\cite{Molnar22}}]\label{prop:Molnar}
\begin{enumerate}[(a)]
  \item  Let $\phi: \bM_n(\mathbb{C})\to \bM_n(\mathbb{C})$ satisfy  \eqref{eq:w2-local+}   and $n\geq 3$.  Then
 $\phi$  is either an algebra automorphism or an algebra anti-automorphism.
    \item  Let $\phi: \bM_n(\mathbb{C})\to \bM_n(\mathbb{C})$ satisfy   \eqref{eq:w2-localx}.  Then
  either $\phi$ or $-\phi$ is an algebra automorphism.
  \item  Let  $\cA$ be a standard
operator algebra over a complex Banach space $E$.  Let $\phi : \cA \to \cB(E)$ satisfy  \eqref{eq:w2-localD}.
Then $\phi$ is a linear derivation.
\end{enumerate}
\end{prop}

In this note,
we  extend above  results  of Moln\'{a}r, among many others,
 about operational $2$-local automorphisms to the setting of standard operator algebras.
We also show that any operational $2$-local derivation of a semi-finite von Neumann algebra
or a unital simple $C^{\ast}$-algebra with a  faithful tracial state is a linear local derivation.


\section{Operational $2$-local automorphisms}

Let $E$ be a (complex)   Banach space.
Denote by  $\cB(E)$    the algebra of  bounded linear operators on $E$,
 and $\cF(E)$  its subalgebra of (continuous)  finite rank operators.
A \emph{standard operator algebra} $\cA$ on  $E$
is a subalgebra of   $\cB(E)$   containing   $\cF(E)$.

A  (linear) Jordan isomorphism $J$ between standard operator algebras
is either an algebra isomorphism or an algebra anti-isomorphism; namely, $J$ assumes either the form
\begin{itemize}
\item $T\mapsto  STS^{-1}$ for  an invertible bounded linear map $S:E_1\to E_2$,  or
\item  $T\mapsto ST'S^{-1}$ for an invertible bounded linear map $S: E_1'\to E_2$,
\end{itemize}
where
$E'$ is the Banach dual space of $E$,
 and $T': E_1'\to E_1'$ is the dual map of $T$.
When the second form holds, both $E_1$ and $E_2$ are reflexive (\cite[Theorem 2.6]{hlw08}).

\begin{prop}\label{prop:2-local+auts}
  Let $\cA$ be  a standard operator algebra $\cA$ on a complex Banach space $E$ of dimension at least three.
  Let  $\phi: \cA \to \cA$ be a map such that the range of $\phi$ contains all rank one idempotents.
  Suppose for any $a,b\in \cA$ there is an algebra automorphism $\theta_{a,b}$ of $\cA$ such that
\begin{align}\label{eq:2-local+}
\phi(a) + \phi(b) = \theta_{a,b}(a+b).
\end{align}
Then
the restriction
$\phi\mid_{\cF(E)}$ is an algebra automorphism or an anti-auto\-mor\-phism
of $\cF(E)$.  If $\phi$ is continuous with respect to the weak operator topology then
$\phi$ is an algebra homomorphism or an anti-homomorphism of $\cA$.
\end{prop}
\begin{proof}
With $b=-a$  in \eqref{eq:2-local+}, we  see that $\phi(-a)=-\phi(a)$ for any $a\in \cA$.
Thus, we have $\phi(0)=0$ and
$$
\phi(a) - \phi(b) = \theta_{a,-b}(a-b) \quad\text{for all $a,b\in \cA$.}
$$
With $b=0$ in \eqref{eq:2-local+}, we see that $\phi$ is a local automorphism of $\cA$.
In particular, $\phi$ sends idempotents to idempotents, and preserves rank.
Let $p,q$ be disjoint  idempotents in $\cA$ of finite rank $m$ and $n$; that is, $pq=qp=0$.
Then \eqref{eq:2-local+} implies that $\phi(p)+\phi(q)$ is an idempotent of rank $m+n$.
Hence, $\phi$ sends exactly disjoint idempotents to disjoint idempotents.
In particular, $\phi$ induces a bijection of the set of all rank-one idempotents
which preserves disjointness in both directions.

Let $p$ be an idempotent of finite rank $n\geq 3$.
 Then $p'=\phi(p)$ is again an idempotent of rank $n$.
Let $a=pap$ in $\cF(E)$.  For any finite rank idempotent $q$ disjoint from $p$, we have
$$
\phi(a) + \phi(q) = \theta_{a,   q}(a + q)
$$
has rank strictly greater than the rank of $a$.
This forces $\phi(a)$  disjoint from $I-\phi(p)= I -p'$.
Indeed,  suppose that $\phi(a)(I-p')$  has rank $1\leq m\leq n$.
We can
write  $\phi(a)(I-p')= \sum_{j=1}^{m} y_j\otimes f_j$ with   linearly independent  vectors $y_1, \ldots, y_m$ in
$E$ and  linearly independent norm one  linear functionals $f_1,\ldots, f_m$ such that   $f_j=0$ on $p'E$ for $j=1,\ldots, m$.
Here,  $v\otimes g$ denotes the rank one operator $u\mapsto g(u)v$.
Let $x_1$ be a unit vector
in $(I-p')E$ such that  $f_1(x_1)=1$ and $f_2(x_1)=\cdots = f_m(x_1)=0$.
 Then
the rank one idempotent $q'= x_1\otimes f_1$ is disjoint from $p'$.
Note that
$$
\phi(a)+q' = \phi(a)(I-p') +q'+ \phi(a)p' = (x_1+y_1)\otimes f_1 + \sum_{j=2}^{m} y_j\otimes f_j + \phi(a)p'
$$
has rank at most the rank of $\phi(a)$.
Let $q$ be a rank one idempotent in $\cF(E)$ such that $\phi(q)=q'$.
Since $\phi(p)\phi(q)=p'q'=0$, we have $pq=0$, and thus $qa=aq = 0$.  We then arrive at a contradiction:
$$
\operatorname{rank}(\phi(a))\geq \operatorname{rank}(\phi(a) + \phi(q)) = \operatorname{rank}(a+q)
= \operatorname{rank}(a)+1 = \operatorname{rank}(\phi(a))+1.
$$
It says $\phi(a)p'=\phi(a)$.  Similarly, we see that $p'\phi(a)=\phi(a)$.
Therefore,
$$
\phi(a) = p'\phi(a)p' \quad\text{for any $a=pap\in \cF(E)$.}
$$
In other words, $\phi$ sends $p\cF(E) p \cong M_n(\mathbb{C})$ into $p'\cF(E)p'\cong M_n(\mathbb{C})$.
Applying Proposition \ref{prop:Molnar}(a), we see that
$\phi$ induces a linear Jordan isomorphism from $p\cF(E)p$ onto $p'\cF(E)p'$.

For any $a, b\in \cF(E)$, let $p$ be an  idempotent of big enough finite rank at least three
 such that $a=pap$ and $b=pap$.  Since $\phi$ induces
 a linear Jordan isomorphism from $p\cF(E)p$ onto $\phi(p)\cF(E)\phi(p)$, we have
 $\phi(\alpha a + \beta b)=\alpha \phi(a) +\beta\phi(b)$ and $\phi(ab + ba)=\phi(a)\phi(b) + \phi(b)\phi(a)$
  for all scalars $\alpha, \beta$.
Since the range of $\phi$ contains
all rank one idempotents and $\phi$ preserves rank, $\phi$ restricts to a linear Jordan automorphism from $\cF(E)$ onto $\cF(E)$.

Finally, since every bounded linear operator on $E$ is a limit
of a net of finite rank operators  in the weak operator topology, the last assertion is established.
\end{proof}

\begin{defn}\label{defn:gen-prod}\label{defn:gen-Jprod}
Fix a finite sequence $(i_1,i_2, \ldots , i_m)$ such that
$\{i_1, i_2, \ldots , i_m\}=\{1, 2, \ldots , k\}$ and $k\geq 2$.
Suppose  there is an index $i_p$ different from all other indices  $i_q$.
Define a \emph{generalized product}
for operators $T_1,\ldots , T_k$ by
$$
T_1*\cdots *T_k\ =\ T_{i_1}\cdots T_{i_m}.
$$
Similarly, we define a \emph{generalized Jordan product}
$$
T_1\circ \cdots \circ T_k
\ =\ T_{i_1} \cdots T_{i_m}\; + \; T_{i_m}\cdots  T_{i_1}.
$$
\end{defn}

The generalized product covers
the usual product $T_1* \cdots * T_k = T_1 \cdots T_k$
and the  Jordan triple product $T_1*T_2 = T_2T_1T_2$.
The generalized Jordan product covers
the usual Jordan product $T_1 \circ T_2 = T_1T_2+T_2T_1$ and the Jordan ternary product
$T_1 \circ T_2 \circ T_3 =T_1T_2T_3 + T_3T_2T_1$.

\begin{thm}[{\cite[Theorem 3.2]{hlw08} and \cite[Theorem 2.2]{hlw10}}]\label{thm:gen-prod}\label{thm:2-local-Jaut}
Let $\mathcal{A}_1, \mathcal{A}_2$ be
standard operator algebras on complex Banach spaces $E_1, E_2$, respectively.
Let $T_1*\cdots *T_k$ and $T_1\circ \cdots \circ T_k$
be a generalized  product and a generalized Jordan product defined as in Definition \ref{defn:gen-prod},
and $\Phi :\mathcal{A}_1\rightarrow\mathcal{A}_2$ is a (not necessarily linear, additive or continuous) map.
 Suppose
 \begin{enumerate}[(a)] \itemsep=5pt
  \item $\Phi(\cA_1)$ contains all finite rank operators on $E_2$ of rank at most two, and the spectrum
\begin{align}\label{eq:gprod}
\sigma(\Phi(A_1)*\cdots *\Phi(A_k)) =\sigma(A_1*\cdots *A_k)
\end{align}
 holds whenever any one of $A_1, \dots, A_k$ has rank at most one; or

\item $\Phi(\cA_1)$ contains all finite rank operators on $E_2$ of rank at most three, and the spectrum
 \begin{align}\label{eq:gJprod}
 \sigma(\Phi(A_1)\circ\cdots \circ\Phi(A_k)) =\sigma(A_1\circ\cdots \circ A_k)
\end{align}
 holds whenever any one of $A_1, \dots, A_k$ has rank at most one.
 \end{enumerate}
Then there exist a scalar $\xi$ with $\xi^{m}=1$ such that $\Phi= \xi J$ for a linear
Jordan homomorphism $J:\cA_1\to \cA_2$.

Moreover, $J$ will be an algebra homomorphism in case (a) unless the generalized product is symmetric in the sense that
$m=2p-1$ and
$$
(i_1, \ldots, i_{p-1}, i_p, i_{p+1}, \ldots, i_{2p-1}) = (i_{2p-1}, \ldots, i_{p+1}, i_p, i_{p-1}, \ldots, i_1).
$$
\end{thm}

\begin{thm}\label{thm:2-local-aut}
Let $\cA$ be  a standard operator algebra  on a complex Banach space $E$.
Let $T_1*\cdots * T_k$ and $T_1\circ\cdots\circ T_k$ be the generalized product and Jordan products of operators in
 $\cA$ defined as in Definition \ref{defn:gen-prod}.
 Let  $\phi: \cA \to \cA$ be a map such that $\Phi(\cA)$ contains all
  finite rank operators on $E$ of rank at most two (resp.\ three).
Suppose
  for any $a_1, \ldots, a_k$ in $\cA$, in which one of them has rank at most one,
   there is an algebra (resp.\ linear Jordan) automorphism $\theta_{a_1,\ldots, a_k}$ of $\cA$ such that
\begin{align}
  \phi(a_1)*\cdots*\phi(a_k) &= \theta_{a_1,\ldots, a_k}(a_1*\cdots * a_k) \label{eq:k-local} \\
  (\text{resp.\quad}    \phi(a_1)\circ\cdots\circ\phi(a_k) &= \theta_{a_1,\ldots, a_k}(a_1\circ\cdots \circ a_k)).  \label{eq:k-localJ}
\end{align}
Then   $\phi$ assumes either the form
$$
A\mapsto \xi SAS^{-1} \quad\text{or}\quad A\mapsto \xi SA'S^{-1}
$$
for  a complex scalar $\xi$ with $\xi^m = 1$, and
an invertible bounded linear operator $S:E\to E$ or $S: E'\to E$.
\end{thm}
\begin{proof}
 Since both algebra isomorphisms and linear Jordan isomorphisms preserve spectrum, condition \eqref{eq:k-local}
  implies condition \eqref{eq:gprod}, and condition \eqref{eq:k-localJ} implies
  condition \eqref{eq:gJprod}.  It then follows from Theorem \ref{thm:gen-prod} the desired assertions.
\end{proof}

\begin{cor}\label{cor:2-local-aut-egs}
  Let $\cA$ be  a standard operator algebra $\cA$ on a complex Banach space.
  Let  $\phi: \cA \to \cA$ be a map such that  the range of $\phi$ contains all
  finite rank operators on $E$ of rank at most two.
  Consider the following operational $2$-local automorphism conditions.
  \begin{enumerate}[(a)]
    \item  $\phi(a)\phi(b) = \theta_{a,b}(ab)$,
    \item $\phi(a)\phi(b)\phi(a) = \theta_{a,b}(aba)$, or
    \item $\phi(a)\phi(b)\phi(a)^2 = \theta_{a,b}(aba^2)$,
    \end{enumerate}
where $\theta_{a,b}$ is an automorphism of $\cA$ for any $a,b\in \cA$, in which one of $a, b$  has rank at most one.

If (a) holds, then either $\phi$ or $-\phi$ is an automorphism.  If (b) holds then one of
$\phi$, $\omega\phi$ or $\omega^2 \phi$ is either an automorphism or an anti-automorphism, where
$\omega$ is a primitive cubic root of unity.  If (c) holds then one of $\pm \phi$ and $\pm \sqrt{-1} \phi$ is
an automorphism.
  \end{cor}
\begin{proof}
  It follows from Theorem \ref{thm:2-local-aut} that $\xi^{-1}\phi$ is an automorphism or an anti-automorphism of $\cA$ for
  some $m$th root $\xi$ of unity with $m=2, 3$ or $4$.  The assertions follow from Theorem \ref{thm:gen-prod}
   since  the products in (a) and (c) are not symmetric.
\end{proof}

   We   note that  $aba$ and its transpose $(aba)^t = a^tb^ta^t$ are similar if $a$ or $b$ has finite rank, when $E$ is a Hilbert space.
  Thus the transpose map $a\mapsto a^t$ satisfies condition (b) in Corollary \ref{cor:2-local-aut-egs}.

It follows from Theorem \ref{thm:2-local-aut} the same way the following results for operational
$2$-local and $3$-local Jordan automorphisms.

\begin{cor}\label{cor:2-local-Jaut-egs}
  Let $\cA$ be  a standard operator algebra $\cA$ on a complex Banach space $E$.
  Let  $\phi: \cA \to \cA$ be a map such that the range of $\phi$ contains all finite rank operators on $E$ of
  rank at most three.
    Consider the following   conditions.
  \begin{enumerate}[(a)]
    \item  $\phi(a)\phi(b) +\phi(b)\phi(a)= \theta_{a,b}(ab+ba)$,
    \item $\phi(a)\phi(b)\phi(c) + \phi(c)\phi(b)\phi(a) = \theta_{a,b,c}(abc+cba)$, or
    \item $\phi(a)\phi(b)\phi(a)^2 + \phi(a)^2\phi(b)\phi(a)= \theta_{a,b}(aba^2 + a^2ba)$,
    \end{enumerate}
where $\theta_{a,b}$ (resp.\ $\theta_{a,b,c}$) is a  linear Jordan automorphism of $\cA$ for any $a,b,c\in \cA$,
in which one of $a, b$ (resp.\ $a,b,c$) has rank at most one.

If (a) holds, then  $\phi$ or $-\phi$ is  a linear Jordan automorphism.  If (b) holds then one of
$\phi$, $\omega\phi$ or $\omega^2 \phi$ is a linear Jordan automorphism, where
$\omega$ is a primitive cubic root of unity.  If (c) holds then one of $\pm \phi$ and $\pm \sqrt{-1} \phi$ is
a linear Jordan automorphism.
  \end{cor}

The assumption that the range of $\phi$ contains enough operators of small ranks is indispensable for $\phi$ being surjective.
Consider, for example, the map $\phi(a)= RaL$ of $\cF(\ell_2)$, where $R$ is the unilateral   shift and $L=R^*$ is the backward
unilateral  shift on the separable Hilbert space $\ell_2$.  Then $\phi$ is a linear $n$-local automorphism of $\cF(\ell_2)$ for
any $n\geq 1$, and thus satisfies all operational $n$-local automorphism conditions.
However, $\phi$ is not an automorphism of $\cF(\ell_2)$.

\begin{thm}[{\cite[Theorem 4.2]{hlw08} and \cite[Theorem 4.1]{hlw10}}]\label{thm:gen-prod-sa}
Let $\cS(H_1), \cS(H_2)$ be the sets of bounded self-adjoint operators on complex Hilbert spaces $H_1, H_2$.
Let $T_1*\cdots *T_k$ and $T_1\circ \cdots \circ T_k$
be a generalized  product and a generalized Jordan product defined as in Definition \ref{defn:gen-prod},
and $\Phi :\mathcal{A}_1\rightarrow\mathcal{A}_2$ is a (not necessarily linear, additive or continuous) map.
 Suppose the range of $\Phi$ contains all finite rank self-adjoint operators on $H_2$ of rank at most two (resp.\ three), and the spectrum
\begin{align*}
\sigma(\Phi(A_1)*\cdots *\Phi(A_k)) &= \sigma(A_1*\cdots *A_k) \\
  (\text{resp.\quad}  \sigma(\Phi(A_1)\circ\cdots \circ\Phi(A_k)) &=\sigma(A_1\circ\cdots \circ A_k))
\end{align*}
 holds whenever any one of the self-adjoint operators $A_1, \dots, A_k$ has rank at most one.
Then there exist a scalar $\xi\in\{-1, 1\}$ with $\xi^{m}=1$ and a surjective linear isometry $U: H_2\to H_1$ such that $\Phi$
assumes either the form $A \mapsto \xi UAU^*$ or $A\mapsto \xi UA^tU^*$.  Here, $A^t$ stands for the transpose of $A$ with respect
to a fixed orthonormal basis of $H_1$.
\end{thm}

Note that we do not need to assume  the generalized product of self-adjoint operators is again self-adjoint
in \cite[Theorem 4.2]{hlw08}.

One can derive some results about operational $2$-local $*$-automorphisms   from Theorem \ref{thm:gen-prod-sa}.
Below is an example.

\begin{cor}\label{cor:2-local-*-aut}
Let $H$ be  a complex Hilbert space.
  Let  $\phi: \cB(H)  \to \cB(H)$ be a map such that the range of $\phi$ contains all
operators in $\cB(H)$ of rank at most two.
Suppose that  for any $a,b\in \cB(H)$
 there is a  $*$-automorphism $\theta_{a,b}$ of $\cA$ such that
\begin{align}\label{eq:local-*aut}
\phi(a)^*\phi(b) = \theta_{a,b}(a^*b).
\end{align}
Then $\phi$ assumes either the form $a\mapsto UaV$  or $a\mapsto  Ua^tV$ for
some unitary operators $U,V$ of $H$.
\end{cor}
\begin{proof}
Putting $a=b=I$
 in  \eqref{eq:local-*aut}, we see that $\phi(I)^*\phi(I)=I$.
Putting $b=a$ in  \eqref{eq:local-*aut}, we see that $\|\phi(a)\|=\|a\|$ for any $a\in \cB(H)$.
For any unit vector $x$ in $H$, let $b\in B(H)$ such that
 $\phi(b)=x\otimes x$ be the self-adjoint rank one operator $y\mapsto \langle y, x\rangle x$.
Together with $a=I$,   condition \eqref{eq:local-*aut} implies
$$
\|\phi(I)^* x\|=\|(\phi(I)^*x)\otimes x\|=\|\phi(I)^*\phi(b)\|= \|b\| = \|x\otimes x\| = 1.
$$
Consequently, the operator $\phi(I)^*$ is an isometry, and thus $\phi(I)\phi(I)^*=I$ as well.
In other words, $\phi(I)$ is a unitary operator on $H$.

Replacing $\phi$ by the map $\phi(I)^*\phi(\cdot)$, we can assume $\phi(I)=I$.
It then follows from \eqref{eq:local-*aut} with $b=I$ that $\phi$ sends every
self-adjoint operator $a$ to a self-adjoint operator $\phi(a)=\theta_{a, I}(a^*)^*$.
Therefore, $\phi$ sends $\cS(H)$ into $\cS(H)$.
Applying Theorem \ref{thm:gen-prod-sa} with the
product $a*b=ab$, we have a unitary operator $W$ of $H$ such that
$\phi(b)= \pm WbW^*$ or $\pm Wb^tW^*$ for every self-adjoint operator $b$ in $\cS(H)$.
Replacing $\phi$ by the map $\pm W^*\phi(\cdot)W$ or $\pm W^*\phi(\cdot)^tW$, we can further assume that
$\phi(b)=b$ for every $b\in \cS(H)$.  For any $x\in H$, it then follows from \eqref{eq:local-*aut}
that the inner product
\begin{align*}
\langle \phi(a)^*x, x\rangle &= \operatorname{trace}(\phi(a)^*(x\otimes  x))
= \operatorname{trace}(\phi(a)^*\phi(x\otimes x)) \\
&= \operatorname{trace}(\theta_{a, x\otimes x}(a^*x\otimes x))
=  \langle a^*x, x\rangle.
\end{align*}
The polar identity of inner products ensures that $\phi(a)=a$ for all $a\in \cB(H)$.

In the original setting, with $U= \pm\phi(I)W$ and $V=W^*$, we arrive at the desired conclusion.
\end{proof}

\section{Operational $2$-local derivations}

Let $\M$ be a semi-finite von Neumann algebra with a normal semi-finite faithful trace $\tau$.
In other words, $\tau : \M_+\to [0, +\infty]$ is a map satisfying that
\begin{enumerate}[1.]
\item $\tau(x+y)=\tau(x) + \tau(y)$ for all $x,y\in \M_+$;
\item $\tau(tx)=t\tau(x)$ for all $x\in \M_+$ and $t\geq 0$;
\item $\tau(x^*x)=\tau(xx^*)$ for all $x\in \M$;
\item $\tau(x^*x)=0$ if and only if $x=0$;
\item $\tau(x) = \sup\, \{\tau(y): 0\leq y\leq x, \tau(y)<+\infty\}$;
\item  $\tau(x)= \sup_i \tau(x_i)$ if $x\in \M_+$ is the $\sigma$-weak limit of an increasing net $x_i$ in $\M_+$.
\end{enumerate}

Set
$$
m_\tau = \{x\in \M: \tau(|x|) < +\infty\},
$$
where $|x|=\sqrt{x^*x}$.
Then $m_{\tau}$ is a $\ast$-algebra and is also a two-sided ideal of $\M$.
The $\sigma$-weak closure of $m_\tau$ is $\M$, and for any
$x\in \M_+$
there exists an increasing net $\{x_\lambda\}_\lambda$ of positive elements in $m_\tau$ with strong (SOT) limit $x$.
Moreover, $\tau$ extends to $m_\tau$ a complex linear map such that $\tau(ab)=\tau(ba)$ for all $a,b\in m_\tau$.
When $\tau(1)=1$, we say that $\tau$ is a normal faithful tracial state; this happens exactly when $\M$ is a finite
von Neumann algebra. See, for example, \cite{Takesaki}.


The following strengthens the  results of \cite{AA} and \cite[Theorem 2.7]{Molnar22}.

\begin{thm}\label{THM1}  Let $\M$ be a semi-finite von Neumann algebra with a normal faithful semi-finite trace $\tau$.
Let $\phi:\M\rightarrow \M$ be a map.
Suppose for any $A, B\in \M$, there is a linear derivation $D_{a, b}:\M\rightarrow \M$ such that $$\phi(a)b+a\phi(b)=D_{a, b}(ab).$$
 Then $\phi$ is a linear derivation.
 \end{thm}
\begin{proof}
Since  $\phi(I)I+I\phi(I)=D_{I, I}(I)=0$, we have $\phi(I)=0$.   Then the condition
  $$
  \phi(a) = \phi(a)I + a\phi(I)= D_{a,I}(a), \quad\forall a\in \M,
  $$
  says that $\phi$ is a local derivation.
%
%
%
Since the linear derivation $D_{a,I}$ is inner, we have $\phi(a)\in m_{\tau}$ whenever
$a\in m_\tau$.  Consequently, $\phi(m_{\tau})\subseteq m_{\tau}$.

For any $x\in \M, y\in m_{\tau}$, there is a linear derivation $D_{x, y}$ on $\M$
such that  $$\phi(x)y+x\phi(y)=D_{x, y}(xy).$$ Since $D_{x, y}$
is an inner derivation, there exists $a\in \M$ such that
$$
D_{x, y}(xy)=axy-xya\in m_\tau.
$$
Consequently,
$$
\tau(D_{xy}(xy))=0=\tau(\phi(x)y+x\phi(y)),
$$
and thus,
$$
\tau(\phi(x)y)=-\tau(x\phi(y)).
$$
For any $t_1, t_2\in \mathbb{C}$, $ a_1, a_2\in \M$, $b\in m_{\tau}$, we have
\begin{align*}
\tau(\phi(t_1a_1+t_2a_2)b)&=-\tau((t_1a_1+t_2a_2)\phi(b)) \\
&=-[t_1\tau(a_1\phi(b))+t_2\tau(a_2\phi(b))]\\
&=\tau([t_1\phi(a_1)+ t_2\phi(a_2)]b).
\end{align*}
This gives
$$
\tau([\phi(t_1a_1+t_2a_2)-t_1\phi(a_1)- t_2\phi(a_2)]b)=0.
$$
Let
$$
w=\phi(t_1a_1+t_2a_2)-t_1\phi(a_1)- t_2\phi(a_2).
$$
We have $\tau(wb)=0$ for all   $b\in m_{\tau}$.
Take an increasing net $\{e_{\alpha}\}_{\alpha}$ of projections
in $m_{\tau}$ such that $e_{\alpha}\uparrow_{\alpha} I$ in $\M$.
Since $e_{\alpha}w^{\ast}\in m_{\tau}$,
we have
$$
\tau(we_{\alpha}w^{\ast})=0 \quad \text{for all $\alpha.$}
$$
Because $we_{\alpha}w^{\ast}\uparrow_{\alpha}ww^{\ast}$ and $\tau$ is normal,
 we have $\tau(ww^{\ast})=0$. Then $w=0$ since $\tau$ is faithful. Therefore,
$$
\phi(t_1a_1+t_2a_2)-t_1\phi(a_1)- t_2\tau(a_2)=0
$$
{for all scalars $t_1,t_2$ and $a_1, a_2\in \M$;}
namely, $\phi$ is a linear map.
Since every linear local derivation of $M$ is a linear derivation, the proof is complete.
\end{proof}


\begin{thm}\label{THM2} Let $\A$ be a unital simple C$^{\ast}$-algebra with a  faithful tracial state
 $\tau$. 
  Let $\phi:\A\rightarrow \A$ be a  map with the property that for any $a, b\in \A$, there is a linear derivation $D_{a, b}:\A\rightarrow \A$ such that $$\phi(a)b+a\phi(b)=D_{a, b}(ab).$$
 Then $\phi$ is a linear derivation.
 \end{thm}
 \begin{proof}
 Recall that every linear derivation of a unital simple $C^*$-algebra is inner \cite{Saki68}.
Argue as in the proof of Theorem \ref{THM1}, we see that $\phi$ is a linear local derivation of $\A$.
Since every linear local derivation of a $C^*$-algebra is a linear derivation  (\cite[Corollary 1]{shulman94}; see also \cite{bresar92, Johnson01}), the proof is complete.
\end{proof}

\section*{Acknowledgement}
The authors would like to express their appreciation  to Lajos Moln\'ar for several helpful comments.
They are also thankful to Jiankui Li for a valuable discussion on  Theorem \ref{THM1}.

Wang is partially supported in part by NSF of China (11871303, 11971463, 11671133) and NSF of Shandong Province (ZR2019MA039, ZR2020MA008).

\end{document}